\documentclass[12pt,twoside]{article}

\date{}

\begin{document}




\centerline{}

\centerline{}

\centerline {\Large{\bf Some Characterizations of Special Curves}}

\centerline{}

\centerline{\Large{\bf in the Euclidean Space $\mathrm{E}^4$}}

\centerline{}

\newcommand{\mvec}[1]{\mbox{\bfseries\itshape #1}}

\centerline{\bf {Melih Turgut}}

\centerline{Department of Mathematics, Buca Educational Faculty, }

\centerline{Dokuz Eyl\"{u}l University, 35160 Buca, Izmir, Turkey.}

\centerline{E-mail: Melih.Turgut@gmail.com}

\centerline{}

\centerline{\bf {Ahmad T. Ali}}

\centerline{Mathematics Department, Faculty of Science, }

\centerline{Al-Azhar University, Nasr City, 11448, Cairo, Egypt.}

\centerline{E-mail: atali71@yahoo.com}

\centerline{}

\newtheorem{Theorem}{\quad Theorem}[section]

\newtheorem{Definition}[Theorem]{\quad Definition}

\newtheorem{Corollary}[Theorem]{\quad Corollary}

\newtheorem{Lemma}[Theorem]{\quad Lemma}

\newtheorem{Example}[Theorem]{\quad Example}

\begin{abstract}
\textbf{\emph{In this work, first, we express some characterizations of helices and ccr curves in the Euclidean 4-space. Thereafter, relations among Frenet-Serret invariants of Bertrand curve of a helix are presented. Moreover, in the same space, some new characterizations of involute of a helix are presented.}}
\end{abstract}


\section{Introduction}
In the local differential geometry, we think of curves as a geometric set of points, or locus. Intuitively, we are thinking of a curve as the path traced out by a particle moving in $\mathrm{E}^4$. So, investigating position vectors of the curves is a classical aim to determine behavior of the particle (curve).

Natural scientists have long held a fascination, sometimes bordering on mystical obsession for helical structures in nature. Helices arise in nanosprings, carbon nanotubes, $\alpha-$helices, DNA double and collagen triple helix, the double helix shape is commonly associated with DNA, since the double helix is structure of DNA \cite{cam}. This fact was published for the first time by Watson and Crick in 1953 \cite{wats}. They constructed a molecular model of DNA in which there were two complementary, antiparallel (side-by-side in opposite directions) strands of the bases guanine, adenine, thymine and cytosine, covalently linked through phosphodiesterase bonds. Each strand forms a helix and two helices are held together through hydrogen bonds, ionic forces, hydrophobic interactions and van der Waals fores forming a double helix, lipid bilayers, bacterial flagella in Salmonella and E. coli, aerial hyphae in actinomycete, bacterial shape in spirochetes, horns, tendrils, vines, screws, springs, helical staircases and sea shells (helico-spiral structures) \cite{chou, cook}.

Helix is one of the most fascinating curves in science and nature. Also we can see the helix curve or helical structures in fractal geometry, for instance hyperhelices \cite{toledo}. In the field of computer aided design and computer graphics, helices can be used for the tool path description, the simulation of kinematic motion or the design of highways, etc. \cite{yang}. From the view of differential geometry, a helix is a geometric curve with non-vanishing constant curvature $\kappa$ and non-vanishing constant torsion $\tau$ \cite{barros}. The helix may be called a {\it circular helix} or {\it W-curve} \cite{ilarslan, mont1}.

Its known that straight line ($\kappa(s)=0$) and circle ($\tau(s)=0$) are degenerate-helix examples \cite{kuhn}. In fact, circular helix is the simplest three-dimensional spirals. One of the most interesting spiral example is $k$-Fibonacci spirals. These curves appear naturally from studying the $k$-Fibonacci numbers $\{F_{k,n}\}^{\infty}_{n=0}$ and the related hyperbolic $k$-Fibonacci function. Fibonacci numbers and the related Golden Mean or Golden section appear very often in theoretical physics and physics of the high energy particles \cite{elnasch1, elnasch2}. Three-dimensional $k$-Fibonacci spirals was studied from a geometric point of view in \cite{falc}.

Indeed, in Euclidean 3-space $\mathrm{E}^3$, a helix is a special case of the {\it general helix}. A curve of constant slope or general helix in Euclidean 3-space is defined by the property that the tangent makes a constant angle with a fixed straight line called the axis of the general helix. A classical result stated by Lancret in 1802 and first proved by de Saint Venant in 1845 (see \cite{struik} for details) says that: {\it A necessary and sufficient condition that a curve be a general helix is that the ratio $\frac{\kappa}{\tau}$ is constant along the curve, where $\kappa$ and $\tau$ denote the curvature and the torsion, respectively}.

The notation of a generalized helix in $\mathrm{E}^3$ can be generalized to higher dimensions in the same definition is proposed but in $\mathrm{E}^n$, i.e., a generalized helix as a curve $\psi:R\rightarrow \mathrm{E}^n$ such that its tangent vector forms a constant angle with a given direction $U$ in $\mathrm{E}^n$ \cite{romero}.

Two curves which, at any point, have a common principal normal vector are called Bertrand curves. The notion of Bertrand curves was discovered by J. Bertrand in 1850. Bertrand curves have been investigated in $\mathrm{E}^n$ and many characterizations are given in \cite{gluck}. Thereafter, by theory of relativity, investigators extend some of classical differential geometry topics to Lorentzian manifolds. For instance, one can see, Bertrand curves in $\mathrm{E}_1^n$ \cite{ekm}, in $\mathrm{E}_1^3$ for null curves \cite{bal}, and in $\mathrm{E}_1^4$ for space-like curves \cite{yiltur}. In the fourth section of this paper, we follow same procedure as in \cite{yiltur}.

In this work, first, we aim to give some new characterizations of helices and ccr curves in terms of recent obtained theorems. Thereafter, we investigate relations among Frenet-Serret invariants of Bertrand curve couples, when one of is helix, in the Euclidean 4-space. Moreover, we observe that Bertrand curve of a helix is also a helix; and cannot be a spherical curve, a general helix and a 3-type slant helix, respectively. We also express some characterizations of involute of a helix. We hope these results will be helpful to mathematicians who are specialized on mathematical modeling.

\section{Preliminaries}
To meet the requirements in the next sections, here, the basic elements of
the theory of curves in the space $\mathrm{E}^4$ are briefly presented (A more
complete elementary treatment can be found in \cite{hacisali}).\newline

Let $\alpha :I\subset R\rightarrow \mathrm{E}^4$ be an arbitrary curve in the
Euclidean space $\mathrm{E}^4$. Recall that the curve $\alpha$ is said to be of
unit speed (or parameterized by arclength function $s$) if $\left\langle
\alpha ^{\prime }(s),\alpha ^{\prime }(s)\right\rangle =1$, where $%
\left\langle .,.\right\rangle $ is the standard scalar (inner) product of $\mathrm{E}^4$ given by
\begin{center}
$\left\langle \xi,\zeta \right\rangle =\xi_{1}\zeta_{1}+\xi_{2}\zeta_{2}+\xi_{3}\zeta_{3}+\xi_{4}\zeta_{4},$
\end{center}
for each $\xi=(\xi_{1},\xi_{2},\xi_{3},\xi_{4})$, $\zeta=(\zeta_{1},\zeta_{2},\zeta_{3},\zeta_{4})\in
\mathrm{E}^4$. In particular, the norm of a vector $\xi\in \mathrm{E}^4$ is given by
\begin{center}
$\left\Vert \xi\right\Vert =\sqrt{\left\langle \xi,\xi \right\rangle }.$
\end{center}
Let $\left\{ T(s),N(s),B(s),E(s)\right\} $ be the moving frame along the
unit speed curve $\alpha $. Then the Frenet-Serret formulas are given by
\cite{gluck,Sh}
\begin{equation}\label{u1}
\left[
\begin{array}{c}
T^{\prime } \\
N^{\prime } \\
B^{\prime } \\
E^{\prime }%
\end{array}%
\right] =\left[
\begin{array}{cccc}
0 & \kappa & 0 & 0 \\
-\kappa & 0 & \tau & 0 \\
0 & -\tau & 0 & \sigma \\
0 & 0 & -\sigma & 0%
\end{array}%
\right] \left[
\begin{array}{c}
T \\
N \\
B \\
E%
\end{array}%
\right] .
\end{equation}
Here $T,N,B$ and $E$ are called, respectively, the tangent, the normal, the
binormal and the trinormal vector fields of the curve. and the functions $%
\kappa (s),\tau (s)$ and $\sigma (s)$ are called, respectively, the first,
the second and the third curvature of a curve in $\mathrm{E}^4$. Also, the
functions $H_{1}=\frac{\kappa }{\tau }$ and $H_{2}=\frac{H_{1}^{\prime }}{%
\sigma }$ are called \textit{Harmonic Curvatures} of the curves in $\mathrm{E}^4$,
where $\kappa \neq 0,\tau \neq 0$ and $\sigma \neq 0.$ Let $\alpha :I\subset
R\rightarrow \mathrm{E}^4$ be a regular curve. If tangent vector field $T$ of $\alpha$ forms a constant angle with unit vector $U$, this curve is called an inclined curve or a general helix in $\mathrm{E}^4$. Recall that, A curve $\psi =\psi (s)$ is called a 3-type slant helix if the trinormal lines of $\alpha$ make a constant angle with a fixed direction in $\mathrm{E}^4$ \cite{turyil}. Recall that, if a regular curve has constant Frenet curvatures ratios, (i.e., $\frac{\tau}{\kappa}$ and $\frac{\sigma}{\tau}$ are constants), then it is called a \textit{ccr-curve} \cite{mont2}. It is worth noting that: the W-curve, in Euclidean 4-space $\mathrm{E}^4$, is a spacial case of a ccr-curve.

Let $\alpha (s)$ and $\alpha ^{\ast }(s)$ be regular curves in $\mathrm{E}^4$. $\alpha(s)$ and $\alpha^{\ast}(s)$ are called Bertrand Curves if for each $s_{0}$, the principal normal vector to $\alpha $ at $s=s_{0}$ is the same as the principal normal vector to $\alpha ^{\ast }(s)$ at $s=s_{0}.$ We say that $\alpha ^{\ast }(s)$ is a Bertrand mate for $\alpha (s)$ if $\alpha (s)$ and $\alpha ^{\ast }(s)$ are Bertrand Curves.

In the same space, in \cite{magden1}, author defined a vector product and gave a
method to establish the Frenet-Serret frame for an arbitrary curve by the
following definition and theorem:
\begin{Definition}
Let $a=(a_{1},a_{2},a_{3},a_{4})$, $b=(b_{1},b_{2},b_{3},b_{4})$ and $%
c=(c_{1},c_{2},c_{3},c_{4})$ be vectors in $\mathrm{E}^4$. The vector product
in $\mathrm{E}^4$ is defined by the determinant
\begin{eqnarray}\label{u2}
a\wedge b\wedge c=\left\vert
\begin{array}{cccc}
e_{1} & e_{2} & e_{3} & e_{4} \\
a_{1} & a_{2} & a_{3} & a_{4} \\
b_{1} & b_{2} & b_{3} & b_{4} \\
c_{1} & c_{2} & c_{3} & c_{4}%
\end{array}%
\right\vert,
\end{eqnarray}
where $e_1,e_2,e_3$ and $e_4$ are mutually orthogonal vectors
(coordinate direction vectors) satisfying equations
$$
e_1\wedge e_2\wedge e_3=e_4,\,\,e_2\wedge e_3\wedge e_4=e_1,\,\,
e_3\wedge e_4\wedge e_1=e_2,\,\,e_4\wedge e_1\wedge e_2=e_3.
$$
\end{Definition}

\begin{Theorem}\label{1}
Let $\alpha=\alpha(t)$ be an arbitrary regular curve in the Euclidean
space $\mathrm{E}^4$ with above Frenet-Serret equations. The Frenet apparatus of
$\alpha$ can be written as follows:
\[
T=\frac{\alpha^{\prime }}{\left\Vert \alpha ^{\prime }\right\Vert },
\]%
\[
N=\frac{\left\Vert \alpha ^{\prime }\right\Vert ^{2}\alpha ^{\prime \prime
}-\left\langle \alpha ^{\prime },\alpha ^{\prime \prime }\right\rangle
\alpha ^{\prime }}{\left\Vert \left\Vert \alpha ^{\prime }\right\Vert
^{2}\alpha ^{\prime \prime }-\left\langle \alpha ^{\prime },\alpha ^{\prime
\prime }\right\rangle \alpha ^{\prime }\right\Vert },
\]%
\[
B=\mu\,E\wedge T\wedge N,
\]%
\[
E=\mu\,\frac{T\wedge N\wedge \alpha ^{\prime \prime \prime }}{\left\Vert
T\wedge N\wedge \alpha ^{\prime \prime \prime }\right\Vert },
\]%
\[
\kappa =\frac{\left\Vert \left\Vert \alpha ^{\prime }\right\Vert ^{2}\alpha
^{\prime \prime }-\left\langle \alpha ^{\prime },\alpha ^{\prime \prime
}\right\rangle \alpha ^{\prime }\right\Vert }{\left\Vert \alpha ^{\prime
}\right\Vert ^{4}}
\]%
\[
\tau =\frac{\left\Vert T\wedge N\wedge \alpha ^{\prime \prime \prime
}\right\Vert \left\Vert \alpha ^{\prime }\right\Vert }{\left\Vert \left\Vert
\alpha ^{\prime }\right\Vert ^{2}\alpha ^{\prime \prime }-\left\langle
\alpha ^{\prime },\alpha ^{\prime \prime }\right\rangle \alpha ^{\prime
}\right\Vert }
\]%
and
\[
\sigma =\frac{\left\langle\alpha ^{(IV)},E\right\rangle}{\left\Vert T\wedge N\wedge \alpha
^{\prime \prime \prime }\right\Vert \left\Vert \alpha ^{\prime }\right\Vert }%
,
\]%
where $\mu$ is taken $-1$ or $+1$ to make $+1$ the determinant of $\left[
T,N,B,E\right] $ matrix.
\end{Theorem}

\section{Some New Results of Helices and Ccr Curves}
In this section we state some related theorems and then we express some important results about helices and ccr curves:
\begin{Theorem}
Let $\alpha =\alpha (s)$ be a regular curve in $\mathrm{E}^4$ parameterized by arclength with curvatures $\kappa ,\tau $ and $\sigma $. Then $\alpha =\alpha (s)$ lies on the hypersphere of center $m$ and radius $r\in \Re^{+}$ in $\mathrm{E}^4$ if and only if
\begin{eqnarray}\label{u3}
\rho ^{2}+\left( \frac{1}{\tau }\frac{d\rho }{ds}\right) ^{2}+\frac{1}{%
\sigma ^{2}}\left[ \rho \tau +\frac{d}{ds}\left( \frac{1}{\tau }\frac{%
d\rho }{ds}\right) \right] ^{2}=r^{2},
\end{eqnarray}
where $\rho =\frac{1}{\kappa }$ \cite{mont2}.
\end{Theorem}
\begin{Theorem}
Let $\alpha =\alpha (s)$ be a regular curve in $\mathrm{E}^4$ parameterized by arclength with curvatures $\kappa$, $\tau$ and $\sigma$. Then $\alpha$ is a generalized helix if and only if
\begin{eqnarray}\label{u4}
H'_2+\sigma\,H_1=0,
\end{eqnarray}
where $H_1=\frac{\kappa}{\tau}$ and $H_2=\frac{1}{\sigma}\,H'_1$ are the Harmonic Curvatures of $\alpha$ \cite{magden2}.
\end{Theorem}
\begin{Theorem}
Let $\alpha =\alpha (s)$ be a regular curve in $\mathrm{E}^4$ parameterized by arclength with curvatures $\kappa$, $\tau$ and $\sigma$. Then $\alpha$ is a type 3-slant helix (its second binormal vector $E$ makes a constant angle with a fixed diretion $U$) if and only if
\begin{eqnarray}\label{u5}
\tilde{H}'_2+\sigma\,\tilde{H}_1=0,
\end{eqnarray}
where $\tilde{H}_1=\frac{\sigma}{\tau}$ and $\tilde{H}_2=\frac{1}{\kappa}\,\tilde{H}'_1$ are (we can called) the Anti-Harmonic Curvatures of $\alpha$ \cite{onder}.
\end{Theorem}
With the aid of the above theorems, one can easily obtain the following important results:
\begin{Theorem}
Let $\alpha =\alpha(s)$ be a helix in $\mathrm{E}^4$ with non-zero curvatures.

{\bf 1.} $\alpha$ can not be a generalized helix

{\bf 2.} $\alpha$ can not be a 3-type slant helix

{\bf 3.} If $\alpha$ lies on the hypersphere $S^3$, the sphere's radius is equal to $\frac{\sqrt{\tau^2+\sigma^2}}{\kappa\,\sigma}$.
\end{Theorem}
\begin{Theorem}
Let $\alpha =\alpha(s)$ be a ccr-curve $\mathrm{E}^4$ with non-zero curvatures $\kappa(s)$, $\tau(s)=a\,\kappa(s)$ and $\sigma(s)=b\,\kappa(s)$. Then

{\bf 1.} $\alpha$ can not be a generalized helix

{\bf 2.} $\alpha$ can not be a 3-type slant helix

{\bf 3.} If $\alpha$ lies on the hypersphere $S^3$, then, if and only if, the following equation is satisfied:
\begin{eqnarray}\label{u6}
f^2+\frac{f^{\prime2}}{4a^2}+\frac{f}{4a^2b^2}(2a^2+f^{\prime\prime})^2=r^2,
\end{eqnarray}
where the function $f=f(s)=\rho^2(s)=\frac{1}{\kappa^2(s)}$.
\end{Theorem}

\section{Bertrand Curve of a Helix}
In this section, we investigate relations among Frenet-Serret invariants of Bertrand curve of a helix in the space $\mathrm{E}^4$.
\begin{Theorem}
Let $\delta=\delta(s)$ be a helix in $\mathrm{E}^4$. Moreover, $\xi$ be Bertrand mate of $\delta$. Frenet-Serret apparatus of $\xi$, $\left\{T_{\xi},N_{\xi},B_{\xi },E_{\xi },\kappa_{\xi},\tau_{\xi},\sigma_{\xi}\right\}$, can be formed by Frenet apparatus of $\delta$ $\left\{T,N,B,E,\kappa ,\tau ,\sigma \right\}$.
\end{Theorem}

\textbf{Proof.}
Let us consider a helix (W-curve, i.e.) $\delta=\delta(s)$.  We may express
\begin{equation}\label{u7}
\xi=\delta+\lambda\,N.
\end{equation}
We know that $\lambda=c=$constant (cf. \cite{hacisali}). By this way, we can write that
$$
\frac{d\xi}{ds_{\xi}}\frac{ds_{\xi}}{ds}=T_{\xi}\frac{ds_{\xi}}{ds}
=(1-\lambda\,\kappa)T+\lambda\,\tau\,B
$$
So, one can have
\begin{equation}\label{u8}
T_{\xi }=\frac{(1-\lambda \kappa )T+\lambda \tau B}{\sqrt{\left( 1-\lambda
\kappa \right) ^{2}+(\lambda \tau )^{2}}},
\end{equation}
and
\begin{equation}\label{u9}
\frac{ds_{\xi }}{ds}=\left\Vert\xi ^{\prime}\right\Vert=\sqrt{\left( 1-\lambda \kappa \right)
^{2}+(\lambda \tau )^{2}}.
\end{equation}
In order to determine relations, we
differentiate:
\begin{equation}\label{u10}
\left.
\begin{array}{c}
\xi^{\prime\prime}=\left[\kappa-\lambda(\kappa^2+\tau^2)\right]
N+(\lambda\,\tau\,\sigma)E,\\
\xi^{\prime\prime\prime}=\kappa\left[\lambda(\kappa^2+\tau^2)
-\kappa\right]T+\tau[\kappa-\lambda(\kappa^2+\tau^2+\sigma^2)]B,\\
\xi^{(IV)}=l_1\,N+l_2\,E
\end{array}
\right.
\end{equation}
where
$$
l_1=\kappa^3(\lambda\,\kappa-1)+\lambda\,\tau^2(2\kappa^2+\tau^2+\sigma^2),
$$
and
$$
l_2=\tau\,\sigma[\kappa-\lambda(\kappa^2+\tau^2+\sigma^2)].
$$
Using the above equations, we can form
$$
\left\Vert\xi^{\prime}\right\Vert^2\,\xi^{\prime\prime}-\left\langle
\xi^{\prime},\xi^{\prime\prime}\right\rangle\,\xi^{\prime}=K^2\left[[\kappa-\lambda(\kappa^2+\tau^2)]
N+(\lambda\,\tau\,\sigma)E\right],
$$
where
$$
K=\sqrt{\left( 1-\lambda \kappa \right)
^{2}+(\lambda \tau )^{2}}.
$$
Therefore, we obtain the principal normal and the first curvature,
respectively,
\begin{equation}\label{u11}
N_{\xi}=\frac{1}{L}\left[[\kappa-\lambda(\kappa^2+\tau^2)]
N+(\lambda\,\tau\,\sigma)E\right],
\end{equation}
and
\begin{equation}\label{u12}
\kappa _{\xi}=\frac{L}{K^2},
\end{equation}
where
$$
L=\sqrt{[\kappa-\lambda(\kappa^2+\tau^2)]^2+(\lambda\,\tau\,\sigma)^2}.
$$
Now, we can compute the vector form $T_{\xi }\wedge N_{\xi }\wedge \xi ^{\prime \prime
\prime }$ as the following:
$$
T_{\xi}\wedge N_{\xi}\wedge \xi^{\prime\prime\prime}=
\begin{array}{c}
\frac{1}{KL}\left\vert
\begin{array}{cccc}
T & N & B & E \\
1-\lambda \kappa & 0 & \lambda \tau & 0 \\
0 & \kappa -\lambda (\kappa ^{2}+\tau ^{2}) & 0 & \lambda \tau \sigma \\
l_{1} & 0 & l_{2} & 0%
\end{array}
\right\vert\\
=-\frac{M}{K\,L}\left[\lambda\,\tau\,\sigma\,N-[\kappa-\lambda(\kappa^2+\tau^2)]E\right]
\end{array}
$$
where
$$
M=\tau\left[\lambda(\kappa^2+\tau^2+\sigma^2)-\kappa(1+\lambda^2\sigma^2)\right].
$$
Since, we have
\begin{equation}\label{u12}
E_{\xi}=-\frac{1}{L}\left[\lambda\,\tau\,\sigma\,N-[\kappa-\lambda(\kappa^2+\tau^2)]E\right].
\end{equation}
By this way, we have the third curvature as follows:
\begin{equation}\label{u13}
\tau _{\xi}=\frac{M}{K^2L}.
\end{equation}
Besides, considering last equation of theorem \ref{1}, one can calculate
\begin{equation}\label{u14}
\sigma _{\xi}=\frac{\kappa\,\sigma}{L}.
\end{equation}
Now, to determine the third vector field of Frenet frame, we write
$$
E_{\xi}\wedge T_{\xi }\wedge N_{\xi }=-\frac{1}{KL^2}\left\vert
\begin{array}{cccc}
T & N & B & E \\
0 & \lambda \tau \sigma & 0 & \lambda(\kappa^2+\tau^2)-\kappa \\
1-\lambda \kappa & 0 & \lambda \tau & 0 \\
0 & \kappa -\lambda (\kappa ^{2}+\tau ^{2}) & 0 & \lambda \tau \sigma
\end{array}
\right\vert,
$$
Since, we obtain:
\begin{equation}\label{u15}
B_{\xi}=-\frac{1}{K}\left[\lambda\,\tau\,T+(1-\lambda\,\kappa)B\right].
\end{equation}
It is worth noting that $\mu=1$.\\

Considering obtained equations, we give:
\begin{Theorem}
Let $\delta= \delta(s)$ be a helix in $\mathrm{E}^4$. Moreover, $\xi$ be Bertrand mate of $\delta$.

{\bf 1.} $\xi$ is also a helix.

{\bf 2.} $\xi$ can not be a generalized helix

{\bf 3.} $\xi$ can not be a 3-type slant helix

{\bf 4.} If $\xi$ lies on the hypersphere $S^3$, then, the sphere's radius is equal to $\frac{\sqrt{\tau_\xi^2+\sigma_\xi^2}}{\kappa_\xi\,\sigma_\xi}=
\frac{\sqrt{\tau^2+(1-\lambda\,\kappa)^2\sigma^2}}{\kappa\,\sigma}$.
\end{Theorem}

\section{Involute-evolute Curve of a Helix}
In this section, first, we correct the computations in the paper \cite{ozy} and then we obtain new results:
\begin{Theorem}
Let $\xi=\xi(s)$ be involute of $\delta$. Let $\delta$ be a helix in $\mathrm{E}^4$. The Frenet apparatus of $\xi$, $\left\{T_{\xi},N_{\xi},B_{\xi },E_{\xi },\kappa_{\xi},\tau_{\xi},\sigma_{\xi}\right\}$, can be formed by Frenet apparatus of $\delta$ $\left\{T,N,B,E,\kappa ,\tau ,\sigma \right\}$ and take the following form.
\begin{equation}\label{u16}
T_\xi=N,\,\,\,\,\,
N_\xi=\frac{-\kappa\,T+\tau\,B}{\sqrt{\kappa^2+\tau^2}},\,\,\,\,\,
B_\xi=-E,\,\,\,\,\,
E_\xi=\frac{\tau\,T+\kappa\,B}{\sqrt{\kappa^2+\tau^2}},
\end{equation}
and
\begin{equation}\label{u17}
\kappa_\xi=\frac{\sqrt{\kappa^2+\tau^2}}{\kappa\,|c-s|},\,\,
\tau_\xi=\frac{\tau\,\sigma}{\kappa\,\sqrt{\kappa^2+\tau^2}\,|c-s|},\,\,
\sigma_\xi=-\frac{\sigma}{\sqrt{\kappa^2+\tau^2}\,|c-s|},
\end{equation}
where
\begin{equation}\label{u18}
\frac{ds_{\xi}}{ds}=\kappa\,|c-s|.
\end{equation}
\end{Theorem}

\textbf{Proof.}
The proof of the above theorem is similar as the proof of the previous theorem.\\

Considering obtained equations, we give:
\begin{Theorem}
Let $\xi$ and $\delta$ be unit speed regular curves in $\mathrm{E}^4$. $\xi$ be involute of $\delta$. Then, the involute

{\bf 1.} $\xi$ cannot be a helix.

{\bf 2.} $\xi$ is a ccr-curve.

{\bf 3.} $\xi$ cannot be a generalized helix

{\bf 4.} $\xi$ cannot be a 3-type slant helix

{\bf 5.} $\xi$ cannot be lies on the hypersphere $S^3$.
\end{Theorem}

\textbf{Proof.}
The proof of points 1, 2, 3 and 4 are obviously. In the following we will proof the point 5:

Integrating the equation (\ref{u18}), we have
$$|c-s|=\sqrt{\frac{2s_\xi}{\kappa}},$$
which leads to
\begin{equation}\label{u19}
\kappa_\xi=\frac{A_1}{\sqrt{s_\xi}},\,\,\,\,\,
\tau_\xi=\frac{A_2}{\sqrt{s_\xi}},\,\,\,\,\,
\sigma_\xi=\frac{A_3}{\sqrt{s_\xi}},
\end{equation}
where
$$
A_1=\sqrt{\frac{\kappa^2+\tau^2}{2\kappa}},\,\,A_2=-\frac{\tau\,\sigma}{2\kappa(\kappa^2+\tau^2)},\,\,
A_3=-\frac{\sigma\sqrt{\kappa}}{\sqrt{2(\kappa^2+\tau^2)}}.
$$
Then if the evolute $\xi$ lies in the hypersphere the equation (\ref{u6}) must be satisfied. Substituting $f=\frac{s_\xi}{A_1^2}$, $\kappa_\xi=\frac{A_1}{\sqrt{s_\xi}}$, $B_1=\frac{\tau_\xi}{\kappa_\xi}$ and $B_2=\frac{\sigma_\xi}{\kappa_\xi}$ in the equation (\ref{u6}), we have
$$
\frac{s_\xi\,(B_1^2+B_2^2)}{A_1^2B_2^2}+\frac{1}{4A_1^2B_1^2}=r^2,
$$
which is contradiction because the radius $r$ of the sphere must be constant and the coefficient of $s_\xi$ can not be equal zero. The proof is completed.

\end{document}